 \documentclass[12pt]{article}     % if you need a4paper
\usepackage{graphicx}
\newtheorem{definition}{Definition}

\newcommand{\EQ}{\begin{equation}\begin{array}{lllllllll}}
\newcommand{\EE}{\end{array}\end{equation}}
\newcommand{\EQQ}{$$\begin{array}{lllllllll}}
\newcommand{\EEE}{\end{array}$$}
\newcommand{\MT}{\left[ \begin{array}{ccccccccc}}
\newcommand{\EM}{\end{array}\right]}
\newcommand{\bean}{\begin{equation}\begin{array}{rcllllllll}}
\newcommand{\eean}{\end{array}\end{equation}}
\newcommand{\bea}{$$\begin{array}{rcllllllll}}
\newcommand{\eea}{\end{array}$$}
\def\ds{\displaystyle}
\def\Fr{\ds \frac}
\def\=={&=&}

\title{\LARGE \bf
Computational Analysis of \\Control Systems \\Using Dynamic Optimization
\thanks{This work was supported in part by U.S. Naval Research Laboratory and Air Force Office of Scientific Research}
}

\author{Wei Kang \thanks{Wei Kang is with Faculty of Applied Mathematics, Naval Postgraduate School, Monterey, CA, USA
        {\tt\small wkang@nps.edu}} 
        and Liang Xu\thanks{Liang Xu is at Naval Research Laboratory, Monterey, CA, USA
        {\tt\small liang.xu@nrlmry.navy.mil}}%
}

\begin{document}

\maketitle
\thispagestyle{empty}
%\pagestyle{empty}

%%%%%%%%%%%%%%%%%%%%%%%%%%%%%%%%%%%%%%%%%%%%%%%%%%%%%%%%%%%%%%%%%%%%%%%%%%%%%%%%
\begin{abstract}

Several concepts on the measure of observability, reachability, and robustness are defined and illustrated for both linear and nonlinear control systems. Defined by using computational dynamic optimization, these concepts are applicable to a wide spectrum of problems. Some questions addressed include the observability based on user-information, the determination of strong observability vs. weak observability, partial observability of complex systems, the computation of $L^2$-gain for nonlinear control systems, and the measure of reachability in the presence of state constraints. Examples on dynamic systems defined by both ordinary and partial differential equations are shown. 

\end{abstract}

%%%%%%%%%%%%%%%%%%%%%%%%%%%%%%%%%%%%%%%%%%%%%%%%%%%%%%%%%%%%%%%%%%%%%%%%%%%%%%%%
\section{Introduction}

Control systems are analyzed and characterized by using fundamental concepts such as observability, reachability, and input-to-output gain \cite{kailath,isidori,zhou}. These concepts have a vast volume of literature. For nonlinear systems, the challenge is to define the concepts so that they are characteristic and fundamental to control systems and, meanwhile, they are practically verifiable. In this paper, the goal is to use dynamic optimization to define quantitative measures of control system properties. Moreover, computational methods of dynamic optimization provide practical tools to numerically implement these concepts in applications.  

In Section \ref{secobs}, the ambiguity in estimation is defined as a measure of observability. This quantity can be numerically computed by solving a dynamic optimization. An example is shown in which some systems observable in traditional sense are not practically observable because of their poor value of ambiguity in estimation. In other words, we can quantitatively tell strongly observable from weakly observable. Another feature of this concept is the capability of taking into account non-sensor information or user knowledge of systems, in addition to the output. For instance, an example is shown in which the system is unobservable under a traditional definition. It turns out that the system is strongly observable if we know the control input has a bounded variation, but without using its accurate upper bound. Moreover, this concept can be used to measure partial observability of complex systems, including the observability of a function of the states and the observability of unknown parameters in a model. 

In Section \ref{secgain}, computational methods for the $L^p$-gain of control systems are introduced. The assumption is that the space of input has finite dimension. Then, the $L^p$-gain can be computed using dynamic optimization. In addition, a method of approximating $L^p$-gain is also introduced, which is based on the eigenvalues of covariance matrices. The methods are exemplified by a nonlinear model of atomic force microscope.

In Section \ref{seccon}, we define the concepts of ambiguity in control and control cost. These definitions take into account the control input as well as systems' constraints. For instance, the concept can be used to quantitatively measure the reachability of nonlinear systems under the constraint that the states must stay in a given region of safety. As an example, heat equation with boundary control is studied.

\section{Observability}
\label{secobs}
Consider a general control system
\EQ
\label{sys}
\dot x= f(t,x,u,\mu), & x \in \Re^{n_x}, & u\in \Re^{n_u} & \mu\in \Re^{n_\mu}\\
y=h(t,x,u,\mu), &y \in \Re^{n_y}\\
z=e(t,x,u,\mu), &z \in \Re^{n_z} \\
(x(\cdot),u(\cdot),\mu) \in {\cal C}
\EE
in which $y$ is the output, $z$ is the variable to be estimated, which is either the state $x$ or a function of $x$ in the case of partial observability for large scale systems. The system state is $x$, $u$ is the control input, $\mu$ is the parameter or model uncertainty. In (\ref{sys}), ${\cal C}$ is a general formulation of constraints. Some examples of constraints include, but not limited to,
\EQQ
E(x(t_0),x(t_f))\leq 0, & \mbox{ end point condition }\\
s(x,u)\leq 0, & \mbox{ state-control constraints }\\
s(x(t_1))=0, & \mbox{ known event at time } t_1\\
\mu_{min} \leq \mu \leq \mu_{max} & \mbox{ model uncertainties }\\
s(x,\mu) =0, & \mbox{ DAE (differential-algebraic equations, $\mu$ is a variable) }\\
\mbox{Variation}(u) \leq V_{max}, & \begin{array}{ll} \mbox{control input with bounded variation}\\ \mbox{(non-sensor information)} \end{array}
\EEE
These constraints represent known information about the system in addition to the measured output $y$. This general form of constraints makes it possible to take into account non-sensor information, or user knowledge about the system, in the estimation process. For instance, some state variables are known to be nonnegative; or a control input has bounded variation; or an event is known to happen at certain moment. All these are valuable information that can be used for the estimation of $z$. The goal of this section is to define a measure for the observability of $z$ using the observation data of $y$ as well as the constraints and the control system model.

\subsection{Definition} 

We assume that variables along trajectories are associated with metrics. For instance, $y=h(t,x(t), u(t),\mu)$, as a function of $t$, has $L^2$ or $L^\infty$ norm; $z=e(t,x(t),u(t),\mu)$ can be measured by its function norm, or by the norm of its initial value $e(t_0,\xi(t_0), u(t_0), \mu)$. A metric used for $z$ is denoted by $|| \cdot ||_Z$; and $||\cdot||_Y$ represents the metric for $y=h(t,x,u,\mu)$. The following definition is applicable to systems with general metrics, including $L^p$ and $L^\infty$. Unless otherwise specified, a norm $||a ||$ for $a\in \Re^k$ is defined by
 $$(a_1^2+\cdots+a_k^2)^{1/2}$$
For any function $h(t)$, $t\in [t_0, t_f]$, its $L^p$-norm is defined by 
$$||h||_{L^p}=\left(\int_{t_0}^{t_f} |h(t)|^p dt\right)^{1/p}$$
The infinity norm is defined by
$$||h||_\infty=\lim_{p\rightarrow \infty} \left(\int_{t_0}^{t_f} |h(t)|^p dt\right)^{1/p}$$
which equals its essential supremum value. In this paper, a triple $(x(t),u(t),\mu)$ represents a trajectory of (1) satisfying the differential equations as well as the constraints. Given a positive number $\epsilon >0$ and a nominal, or true, trajectory $(x(t), u(t),\mu)$. Define 
\bean
\label{estimation}
{\cal E} =\left\{ (\hat x(t),\hat u(t), \hat \mu) | \; ||h(t,\hat x(t),\hat u(t),\hat \mu)-h(t,x(t),u(t),\mu)||_Y \leq \epsilon\right\}
\eean
The number $\epsilon$ is used as an output error bound. If $h(\hat x(t),\hat u(t),\hat \mu)$ stays in the  $\epsilon$ neighborhood of the nominal output $h(t,x(t),u(t),\mu)$, then we consider the trajectory $(\hat x(t),\hat u(t),\hat \mu)$ not distinguishable from the nominal one using output measurement. In this case, any trajectory in $\cal E$ can be picked by an estimation algorithm as an approximation of the true trajectory $(x(t),u(t),\mu)$. For this reason, a trajectory in $\cal E$ is called an estimation of $(x(t),u(t),\mu)$. Similarly, $\hat z=e(t,\hat x,\hat u,\hat \mu)$ is an estimation of $z=e(t,x,u,\mu)$. 

\begin{definition} \label{defobs}
Given a trajectory $(x(t),u(t),\mu)$, $t \in [t_0,t_1]$. Let $\epsilon>0$ be the output error bound. Then the number $\rho_o(\epsilon)$ is defined as follows
\EQ
\label{rho}
\rho_o(\epsilon) = \ds\max_{(\hat x(t),\hat u(t), \hat \mu)}||e(t,\hat x(t),\hat u(t), \hat \mu)-e(t,x(t),u(t),\mu)||_Z\\
\hspace{-0.3in} \mbox{subject to}\\
||h(\hat x(t),\hat u(t), \hat \mu) -h(t,x(t),u(t),\mu)||_Y \leq \epsilon\\
\dot {\hat x}= f(t,\hat x,\hat u, \hat \mu), \\
(\hat x(\cdot), \hat u(\cdot), \hat \mu) \in {\cal C}
\EE
The number $\rho_o(\epsilon)$ is called the ambiguity in the estimation of $z$ along the trajectory $(x(t),u(t),\mu)$. 

Let $U$ be an open set in $(x,u,\mu)$-space and $[t_0,t_1]$ be a time interval. Then the largest value of ambiguity along all trajectories in $U$ is called the ambiguity in the estimation of $z$ in the region $U$. 

\end{definition}

\noindent\underline{\it Remarks}

1. The ratio $\rho_o(\epsilon)/\epsilon$ measures the sensitivity of estimation to the noise in $y$. A small sensitivity value implies strong observability of $z$ in the presence of sensor noise. 

2.  The ratio $\rho_o(\epsilon)/\epsilon$ is closely related to the observability gramian. Consider a linear system
$$\dot x=Ax, \;\; y=Cx$$
Suppose $z=x$ and suppose $||\cdot||_Y$ is the $L^2$-norm. Let $P$ be the observability gramian \cite{kailath}, \cite{krener}, then 
\EQ
\label{obgramian}
||y||_Y^2=x_0^T P x_0
\EE
Given $||y||_Y=\epsilon$, $\rho_o$ equals the maximum value of $||x_0||$ satisfying (\ref{obgramian}). In this case, $x_0$ is an eigenvector of $P$ associated to the smallest eigenvalue $\lambda_{\min}$. We have
\EQ
\label{ob2}
\epsilon^2=\lambda_{min} \rho_o^2
\EE
Therefore, the ratio $\rho_o(\epsilon)^2/\epsilon^2$ equals the reciprocal of the smallest eigenvalue of the observability gramian. For nonlinear systems, one can use observability gramian to approximate $\rho_o(\epsilon)/\epsilon$. An advantage of this approach is that the gramian can be computed empirically without solving the optimization problem (\ref{rho}). Details on empirical computational algorithms for the gramian of nonlinear systems can be found in \cite{krener} and \cite{lall}.   

3. Given a fixed number $\rho_o>0$, from (\ref{obgramian}) the least sensitive initial state defined in (\ref{ob2}), i.e. the eigenvector of length $\rho_o$ associated to $\lambda_{min}$, can be found by the following optimization
$$\arg\min_{||x_0||_X=\rho_o} ||y||_Y=x_0^TPx_0$$
Extending this idea to nonlinear systems, the least observable direction in initial states can be defined as follows: given $\rho_o>0$, let $\epsilon$ be the minimum value from the following problem of minimization
\EQQ
\epsilon=\ds\min_{(\hat x(t),\hat u(t), \hat \mu)}||h(t,\hat x(t),\hat u(t), \hat \mu)-h(t,x(t),u(t),\mu)||_Z\\
\hspace{-0.3in} \mbox{subject to}\\
||\hat x_0-x_0||_X =\rho_o\\
\dot {\hat x}= f(t,\hat x,\hat u, \hat \mu), \\
(\hat x(\cdot), \hat u(\cdot), \hat \mu) \in {\cal C}
\EEE
The resulting initial $x_0$ represents the least observable state on the sphere $||\hat x_0-x_0||_X =\rho_o$ and the ratio $\rho_o/\epsilon$ measures the unobservability of initial states. This definition is a reverse process of Definition \ref{defobs}. However, in some cases it is easier to handle the constraint $||\hat x_0-x_0||_X =\rho_o$ than the inequality of $h(t,x,u,\mu)$ in (\ref{rho}).
 
4. The metric for output in Definition \ref{defobs} can be a vector valued function which is bounded by a vector $\epsilon$. This flexibility is useful for systems using different types of sensors with different accuracy. 

5. Definition \ref{defobs} is independent of estimation methods. It characterizes a fundamental attribute of the system itself, not the accuracy of a specific estimation method. In the following, we compare Definition \ref{defobs} to traditional definitions of observability. It is shown that simple linear systems observable in the traditional sense might be weakly observable or practically unobservable under Definition \ref{defobs}; and, on the other hand, some systems not observable under traditional definitions are practically observable with a small ambiguity in estimation. 
$\diamond$\\

\subsection{Computational dynamic optimization}

The problem defined by (\ref{rho}) is a dynamic optimization. To apply Definition \ref{defobs}, this problem must be solved. Obviously, an analytic solution to (\ref{rho}) is very difficult to derive, if not impossible, especially in the case of nonlinear systems. However, there exist numerical approaches that can be used to find its approximate solution. For instance, various numerical methods are discussed in detail in
\cite{Ho}, \cite{bryson}, and \cite{polak:book}. Surveys on numerical methods for solving nonlinear optimal control problems can be found in \cite{betts:survey,polak73}. The computational algorithm used in this paper is from a family of approaches called direct method
\cite{betts:book, EKR1995, H2000, FR, KGR}. The essential idea of this method is to
discretize the optimal control problem and then solve the resulting finite-dimensional optimization problem. The
simplicity of direct methods makes it an ideal tool for a wide variety
of applications of dynamic optimization with constraints, including (\ref{rho}) in Definition \ref{defobs}.

More specifically, all simulations in this paper use a pseudospectral optimal control method. In this approach, a set of nodes is selected using either the zeros or the critical points of orthogonal polynomials, in our case the Legendre-Guass-Lobatto nodes. Then, the problem of dynamic optimization is discretized at the nodes to result in a nonlinear programming, which is solved using sequential quadratic programming. Details are referred to \cite{EKR1995,FR,KGR}. In some of the following examples, dynamic optimizations are solved using the software package DIDO \cite{dido}. 

A frustration in nonlinear programming is the difficulty of finding global optimal solutions within a given domain. This is no exception in this paper. In all examples, a variety of initial guesses are used to gain a comfortable level of confidence that the result is not a local optimal solution. However, for all examples of nonlinear systems in this paper, the computation cannot guarantee global maximum value for (\ref{rho}). Nevertheless, in the case that a result is not the global maximum value, it still provides a lower bound of the ambiguity value $\rho_o$.   

\subsection{Examples}
For the rest of this section, we illustrate Definition \ref{defobs} using several examples. In the first example, it is shown that the traditional concept of observability is ineffective for systems with large dimensions. It justifies the necessity of a quantitative definition of observability, such as Definition \ref{defobs}. \\
\\
{\it Example}. Consider the following linear system
\bean
\label{examplinear}
&&\dot x_1=x_2\\
&&\dot x_2=x_3\\
&&\vdots\\
&&\dot x_n=-\ds\sum_{i=1}^n \left(\begin{array}{ccc} n\\ i-1\end{array}\right) x_i\\
&& y=x_1
\eean  
Under a traditional definition of observability, this system is perfectly observable for any choice of $n$, i.e. given an output history $y=x_1(t)$, it corresponds to a unique initial state $x_0$. However, if Definition \ref{defobs} is applied to measure the observability, it is a completely different story when the dimension is high. 

Suppose the goal is to estimate $x_0$. We can define $z=x(t)$. Definition \ref{defobs} is applicable with arbitrary metrics. To measure the observability of the initial state, we can use the norm of $x(0)$ as the metric for $z$, i.e.
$$||z(t)||_Z=||x(0)||$$
For this example, the output accuracy is measured by $L^\infty$-norm, 
$$|| y(t) ||_Y=\ds\max_{t\in [t_0,t_f]} |y(t)|$$
Let us assume that the true initial state is 
$$x_0=\MT 0&0&\cdots&1\EM^T$$ 
Let the output error bound be small, $\epsilon=10^{-6}$. So, we assume very accurate observation data. The time interval is $[0,15]$. Problem (\ref{rho}) has the following form
\EQ
\label{problemlinear}
\rho_o(\epsilon) = \ds\max_{\hat x}||\hat x(0)-x(0)||\\
\hspace{-0.3in} \mbox{subject to}\\
||\hat x_1(t)-x_1(t)||_Y \leq \epsilon\\
\dot {\hat x}= f(\hat x)
\EE
It is solved to compute $\rho_o(\epsilon)$. Table \ref{tab:ErrorTolerance} lists the result for $n=2,3,\cdots, 9$.
\begin{table}[ht]
	\centering
		\begin{tabular}{|c|c|c|c|c|}
			\hline
			n&2&3&4&5\\ \hline
			$\rho_o(\epsilon)$&$4.70\times 10^{-6}$&$2.67\times 10^{-5}$&$1.53\times 10^{-4}$&$8.89\times 10^{-4}$\\ \hline	$\epsilon$&$10^{-6}$&$10^{-6}$&$10^{-6}$&$10^{-6}$\\ \hline
			n&6&7&8&9\\ \hline
			$\rho_o(\epsilon)$&$5.20\times 10^{-3}$&$3.01\times 10^{-2}$&$1.75\times 10^{-1}$&$1.02$\\ \hline	$\epsilon$&$10^{-6}$&$10^{-6}$&$10^{-6}$&$10^{-6}$\\ \hline
		\end{tabular}
	\caption{Observability}
	\label{tab:ErrorTolerance}
\end{table}

\begin{figure}
	\begin{center}
		\includegraphics[width=3.0in]{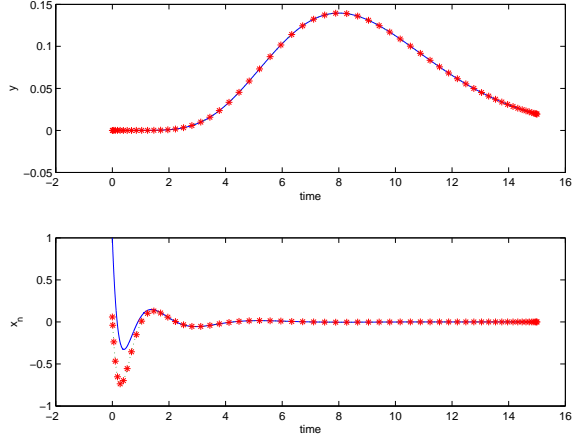} 
		\caption{Estimation error ($n=9$)}
		\label{Estimation error}
		\end{center}
\end{figure}

From the table, when $n=2$ the   ambiguity in the estimation of $x_0$ is as small as $4.70\times 10^{-6}$. So the system is strongly observable. Equivalently, if the observation data has absolute error less than $\epsilon$, then the worst possible estimation of $x_0$ has an error at the scale of $10^{-6}$. This conclusion agrees with the traditional theory of observability. However, when the dimension is increased, the observability ambiguity increases too; thus the system becomes less observable. At $n=8$, the observability ambiguity is as big as $0.175$, or the worst error of estimation is $17.5\%$ relative to the true $x_0$. When $n=9$, the observability ambiguity is $1.02$. In this case, the worse relative error in estimation is more than $100\%$! Thus, the system is practically unobservable, although it is perfectly observable under a traditional definition. Figure \ref{Estimation error} shows why this system is practically unobservable. The continuous curves represent the true trajectory and its output for $n=9$; the dotted curves are the estimation. The outputs of both trajectories agree to each other very well (Figure on top), but the initial states (only $x_9$ is plotted) are significantly different.    

As shown in Figure \ref{Estimation error}, while the estimation of the initial state is inaccurate, the estimation is very close to the true value at the final time $t_f$. To see the observability of the final state, let us use a different metric for $z$,
\EQ\label{newmetric} ||z(t)||_Z=||x(t_f)||\EE
If we consider $t=t_f$ as the current time moment, then this metric is used to measure the detectability of the current system state, rather than the observability of the initial value $x_0$. To compute the ambiguity under the new metric, we solve the problem defined in (\ref{problemlinear}) except that the cost function is replaced by the metric (\ref{newmetric}). For the case of $t_f=10$ and $\epsilon=10^{-6}$, the ambiguity in the estimation of $x(t_f)$ equals $2.7328\times 10^{-6}$. Therefore, the system is accurately detectable.

To summarize, this example shows a set of linear systems that are observable under conventional definition. However, as the dimension is increased, the systems become practically unobservable in the sense that an output trajectory cannot accurately determine the state trajectory. Meanwhile, the detectability of the system is not changed with the dimension. Definition \ref{defobs} is used here to treat both observability and detectability in the same framework, quantitatively.
$\diamondsuit$

\begin{figure}
	\begin{center}
		\includegraphics[width=3.0in]{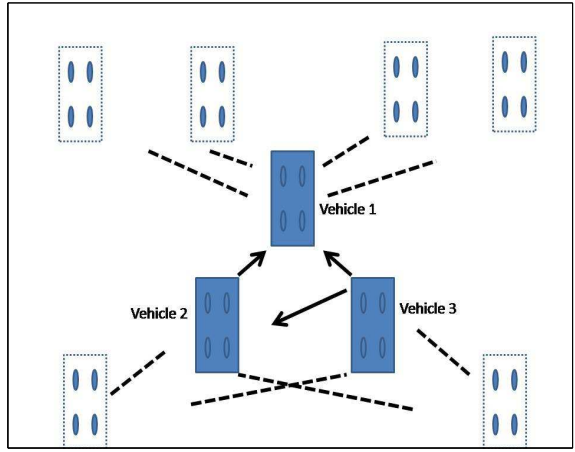} 
		\caption{Cooperative networked system}
		\label{fig:networkedsystem}
		\end{center}
\end{figure}

A concept is useful only if it is verifiable for a wide spectrum of systems and applications. An advantage of Definition \ref{defobs} is that the dynamic optimization (\ref{rho}) can be numerically solved for various types of applications. In the following, we illustrate the usefulness of the ambiguity in estimation using two examples. One is a networked cooperative control system; and the other one is parameter identification for nonlinear systems. 
\\
\\
{\it Example} (Partial observability of cooperative and networked systems). In this example, it is shown that an unobservable system under traditional control theory can be practically observable by employing user knowledge about the system, such as an approximate upper bound of control input. Consider a networked control system showing in Figure \ref{fig:networkedsystem}. Suppose it consists of an unknown number of vehicles. Due to the large number of subsystems, it could be either impossible or unnecessary to process or collect all information about the entire system. A practical approach is to find partial observability with local sensor information only. In this example, we assume that the cooperative relationships in the system is unknown except that we know Vehicle 2 follows Vehicle 1; Vehicle 3 follows both Vehicle 1 and 2 as shown by the arrows in Figure \ref{fig:networkedsystem}. The dashed lines in the figure represent unknown cooperative relationships. The question to be answered is the observability of Vehicle 1 if the locations of Vehicle 2 and 3 can be measured. 

Suppose each vehicle can be treated as a point mass with a linear dynamics
\EQQ
\dot x_{i1}=x_{i2} & \dot y_{i1}=y_{i2}\\
\dot x_{i2}=u_i & \dot y_{i2}=v_i
\EEE
Assume that the control input of vehicles 2 and 3 are defined as follows
\EQQ
u_2=a_1(x_{21}-x_{11}-d_1)+a_2(x_{22}-x_{12})\\
u_3=b_1(x_{31}-\Fr{x_{11}+x_{21}}{2}-d_2)+b_2(x_{32}-\Fr{x_{12}+x_{22}}{2})
\EEE
where $d_i$ is the distance of separation. The control in the $y$-direction is the same. So, Vehicle 2 follows Vehicle 1, Vehicle 3 follows the average position of Vehicle 1 and 2. Suppose we can measure the positions of Vehicles 2 and 3. 
\EQ
\label{mobileoutput}
\mbox{output}=\MT x_{21}&y_{21} & x_{31}&y_{31}\EM^T
\EE
The question to be answered is the observability of the location and velocity of Vehicle 1, i.e. $x_{11}$, $y_{11}$, $x_{12}$, $y_{12}$. We would like to emphasize that the control input of Vehicle 1 is unknown because its input is determined by its cooperative relationships with other vehicles or agents, which is not given. Therefore, in traditional control theory, Vehicle 1 is unobservable. 

To make Vehicle 1 practically observable with limited local measurement, we assume that the input of Vehicle 1 has bounded variation with an upper bound $V_{max}$. This is to say that the vehicles are not supposed to make high frequency zigzag movement, or the control does not have chattering phenomenon. However, discontinuity in control, such as bang-bang, is allowed. In the following, we measure the observability of Vehicle 1 along a trajectory defined by the parameters in Table \ref{NominalTrajectory}.
\begin{table}
	\centering
		\begin{tabular}{|c|c|c|c|c|c|c|c|c|c|c|}
			\hline
			$t_0$&$t_f$&$d_1$&$d_2$&$a_1$&$a_2$&$b_1$&$b_2$&$(x_{11}^0, x_{12}^0)$&$(x_{21}^0, x_{22}^0)$&$(x_{31}^0, x_{32}^0)$\\
			\hline
			$0$&$20$&$-2$&$-2$&$-1$&$-2$&$-3$&$-7$&($0,4$)&($d_1,4$)&($d_2,4$)\\
			\hline
		\end{tabular}
	\caption{Parameters of Nominal Trajectory}
	\label{NominalTrajectory}
\end{table}
The control input of the nominal trajectory is 
$$u_1=\sin\Fr{(t_f-t_0)t}{\pi}$$
which is unknown to the observer. To measure the ambiguity in the estimation of vehicles, we assume the output error bound of (\ref{mobileoutput}) is $\epsilon=10^{-2}$. For the unknown $u_1$, we assume a bounded variation of less than or equal to $V_{max}=3.0$, which is $50\%$ higher than the true variation. The metric for each output variable is the $L^\infty$-norm. The metric for the location and velocity of Vehicle 1 is the $L^2$-norm. The ambiguity in the estimation of each state variable is computed by solving a problem of dynamic optimization. Using the estimation of $x_{11}$ as an example, we have 
\EQ
\label{problemmobile}
\rho_o(\epsilon) = \ds\max_{(\hat x, \hat u_1)}||\hat x_{11}(t)-x_{11}(t)||_{L^2}\\
\hspace{-0.3in} \mbox{subject to}\\
||\hat x_{21}(t)-x_{21}(t)||_{L^\infty} \leq \epsilon_1\\
||\hat x_{31}(t)-x_{31}(t)||_{L^\infty} \leq \epsilon_2\\
\dot {\hat x}= f(\hat x)\\
V(\hat u_1) \leq V_{max}
\EE
where $V(\hat u_1)$ is the total variation. In computation, this constraint is discretized at a set of node points $t_0 < t_1 < \cdots < t_N=f_f$ so that
$$\ds\sum_{k=1}^N |\hat u_1(t_i)-\hat u_1(t_{i-1})| \leq V_{max}$$
An interesting point in the formulation (\ref{problemmobile}) is that the outputs for the two vehicles, i.e. $\hat x_{21}$ and $\hat x_{31}$, have different error bounds, $\epsilon_1$ and $\epsilon_2$. The metric for the outputs is a vector valued function. This flexibility of using different $\epsilon$ value for multiple outputs is advantageous for systems with multiple sensors of different qualities that measure various states with different accuracy. 

\begin{figure}
	\begin{center}
		\includegraphics[width=3.0in]{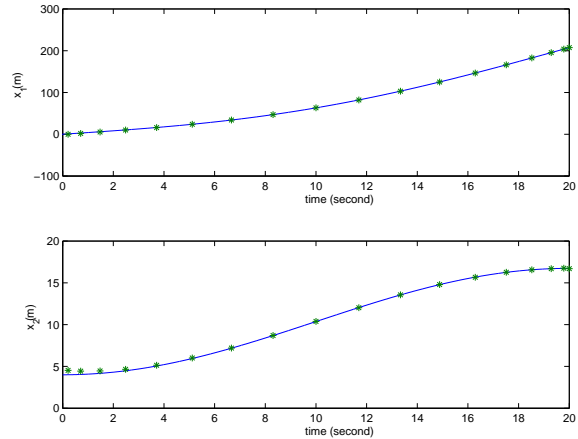} 
		\caption{The worst estimation of the position and velocity of Vehicle 1}
			\label{fig:flow estimation123}
		\end{center}
\end{figure}

The computed result is shown in Table \ref{Observability}. The small relative ambiguity value shows that the location and velocity of Vehicle 1 are practically observable given the measurement of the positions of Vehicle 2 and 3, without using any information about the rest of the networked system and without knowing the input of Vehicle 1. The worst estimation of $x_{11}$ and $x_{12}$ is shown in Figure \ref{fig:flow estimation123}, which has good accuracy.
 $\diamondsuit$
 
\begin{table}
	\centering
		\begin{tabular}{|c|c|c|c|c|c|}
			\hline
			$V_{max}$&$\epsilon$&$\rho_{x_{11}}$&$\rho_{x_{11}}/||x_{11}||_{L_2}$&$\rho_{x_{12}}$&$\rho_{x_{12}}/||x_{12}||_{L_2}$ \\
			\hline
			$3$&$10^{-2}$&$1.2257$&$2.8\times 10^{-3}$&$0.5901$&$1.16\times 10^{-2}$\\
			\hline
		\end{tabular}
	\caption{Observability of Vehicle 1}
	\label{Observability}
\end{table}

\begin{figure}
	\begin{center}
		\includegraphics[width=2.5in]{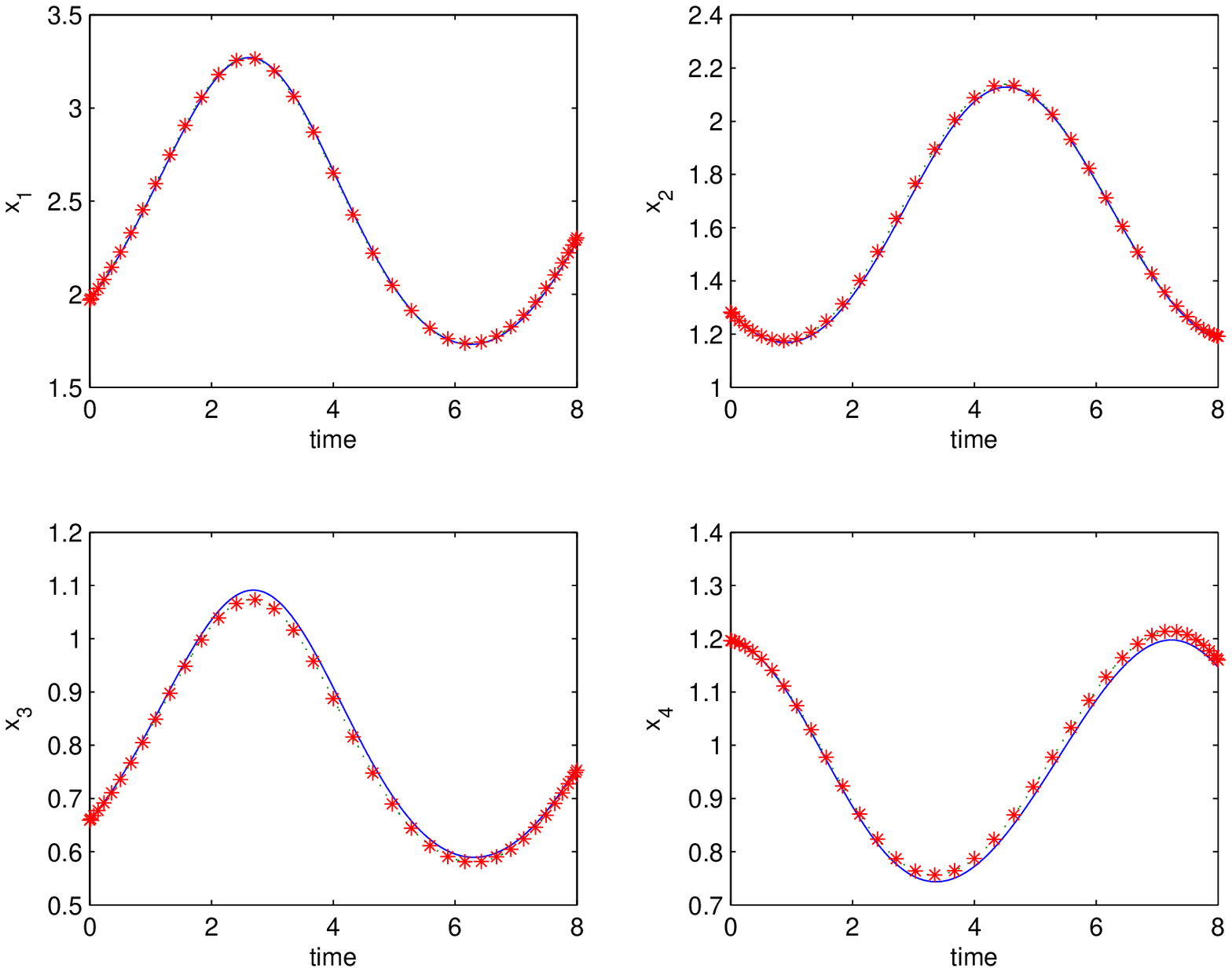} \includegraphics[width=2.5in]{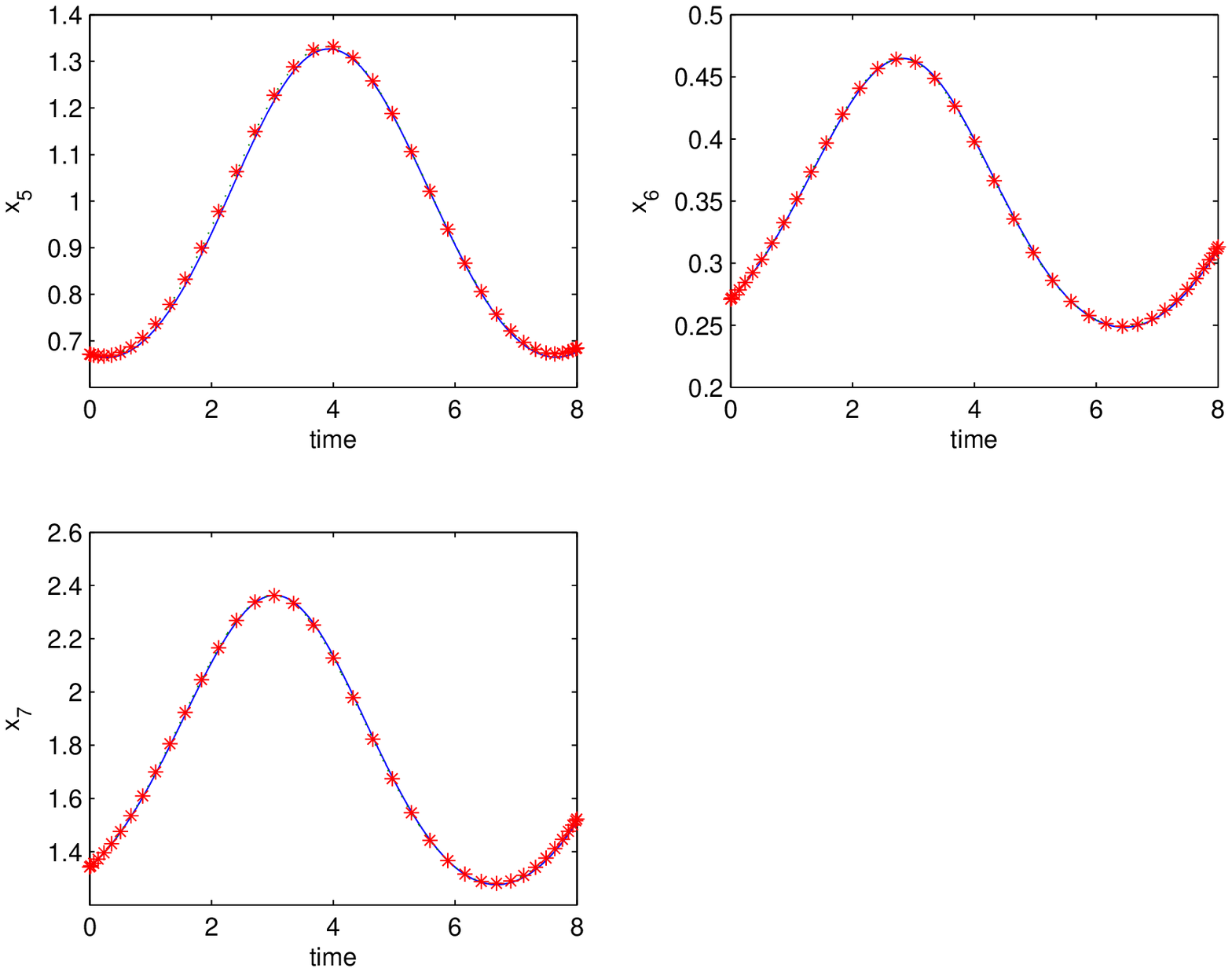}
		\caption{The trajectory of worst estimation: curve - true trajectory; star - estimation)}
		\label{biochemical}
		\end{center}
\end{figure}

In the following example, the concept of ambiguity in estimation is applied to the Laub-Loomis model \cite{laub} with unknown parameters, a nonlinear system of oscillating biochemical network. \\
\\
{\it Example} (Parameter identification) In the study of biochemical networks, it was proposed that interacting proteins could account for the spontaneous oscillations in adenylyl cyclase activity that was observed in homogeneous populations of dictyostelium cells. While a set of terminologies such as $3', 5'$-cycle monophosphate (cAMP) and adenylate cyclase (ACA) are involved in the problem, we focus on the state space in which a set of seven nonlinear differential equations are used as the model \cite{laub}.
\bean
\dot x_1\==k_1x_7-k_2x_1x_2\\
\dot x_2\==k_3x_5-k_4x_2\\
\dot x_3\==k_5x_7-k_6x_2x_3\\
\dot x_4\==k_7-k_8x_3x_4\\
\dot x_5\==k_9x_1-k_{10}x_4x_5\\
\dot x_6\==k_{11}x_1-k_{12}x_6\\
\dot x_7\==k_{13}x_6-k_{14}x_7\\
\eean
In a robustness study \cite{kim}, it was shown that a small variation in the model parameters can effectively destroy the required oscillatory dynamics. As a related question, it becomes interesting to investigate the possibility of estimating the parameters in the system, $k_1$, $k_2$, $\cdots$, $k_7$. To exemplify the idea, we assume that $x_1$, the value of CAC, is measurable, i.e.
$$y=x_1$$
We also assume that the initial states in experimentation is known. Suppose the unknown parameters are $k_1$, $k_6$, and $k_{10}$; and suppose the other parameters are known. The goal is to use the measured data of $y$ to estimate the unknown parameters , i.e.
$$z=\MT k_1 & k_6&k_{10}\EM$$
Using Definition \ref{defobs} we can quantitatively determine the observability of the unknown parameters. Along a nominal trajectory $(x^\ast(t),k^\ast)$, the ambiguity can be computed by solving the following special form of (\ref{rho}). 
\EQQ
\rho_o^2=\displaystyle{\max_{( x,k_1,k_6,k_{10})}} (k_1-k^\ast_1)^2+(k_6-k^\ast_6)^2+(k_{10}-k^\ast_{10})^2\\
\hspace{-0.3in} \mbox{subject to}\\
||x_1(t)-x^\ast_1(t)||_{L^2}^2 \leq \epsilon^2\\
\dot { x}= f(t, x, k_1,k_6,k_{10}), \mbox{ other parameters equal nominal value}\\
x(t_0)=x^\ast(t_0)\\
\EEE
In the simulation, a nominal trajectory is generated using the following parameter value and initial condition 
\EQQ
k_1=2.0,\;\; k_2=0.9,\;\;k_3=2.5,\;\;k_4=1.5,\;\;k_5=0.6,\;\;k_6=0.8,\;\;k_7=1.0,\\
k_8=1.3,\;\;k_9=0.3,\;\;k_{10}=0.8,\;\;k_{11}=0.7,\;\;k_{12}=4.9,\;\;k_{13}=23.0,\;\;k_{14}=4.5,\\
x(0)=\MT 1.9675
    &1.2822
    &0.6594
    &1.1967
    &0.6712
    &0.2711
    &1.3428  \EM
\EEE
The output error bound is being set at $\epsilon = 10^{-2}$. The computation reveals that the ambiguity in the estimation of $z=\MT k_1 & k_6&k_{10}\EM$ is
$$\rho_o=2.38\times 10^{-2}$$
Given the nominal value of the parameters, the relative ambiguity in estimation is about $1\%$. So, the parameters are strongly observable. In fact, the worst estimation of the parameters is 
$$ k_1=2.0150, \;\;  k_6=0.8082, \;\; k_{10}=0.7836.$$
The trajectory generated by the worst parameter estimation is shown in Figure \ref{biochemical}. $\diamond$

\section{Input-to-output gain}
\label{secgain}
$L^p$-gain is a tool of analysis widely used by control engineers to quantitatively measure the sensitivity and robustness of systems. 
Consider
\bean
\label{sysgain}
\dot x=f(t,x,w,\mu)\\
z=e(t,x,w,\mu)
\eean
where $x\in \Re^{n_x}$ is the state variable, $w\in \Re^{n_w}$ is the input that represents the disturbance, $\mu \in \Re^{n_\mu}$ is the system uncertainty or a parameter, $z\in \Re^{n_z}$ is the performance. 
In the following, the $L^p$-norm of a vector valued function is denoted by $||\cdot ||_{L^p}$, for instance
$$|| z(t) ||_{L^p}=\ds\left( \int_{t_0}^{t_1} \sum_{i=1}^{n_z} |z_i(t)|^p dt \right)^{1/p}$$
Given a fixed time interval $[t_0, t_1]$ and $\sigma >0$. Suppose the input $w$ is a function in $L^p$ space for some $1 \leq p \leq \infty$ such that $w(t)$ is bounded by $\sigma$, i.e.
$$ ||w(t)||_{L^p} \leq \sigma$$
Let $x^\ast(t)$ be a nominal trajectory with $x^\ast(t_0)=x_0$, $w^\ast (t)=0$. Fix the initial value $x(t_0)=x_0$. Suppose the system uncertainty is bounded, $\mu_{min} \leq \mu \leq \mu_{max}$. Then the $L^p$-gain from $w$ to $z$ along $x^\ast(t)$ is defined as follows
\EQ
\label{gain}
\gamma (\sigma)=\ds\max_{\begin{array}{c} ||w||_{L^p}\leq \sigma,\\ \mu_{min} \leq \mu \leq \mu_{max}\end{array}} \Fr{||e(t,x,w,\mu)-e(t,x^\ast,0, \mu)||_{L^p}}{\sigma}
\EE

\noindent {\it Remark 6}.
Without parameter, the maximum value of $||e(t,x,w)-e(t,x^\ast,0)||_{L^p}$ is the ambiguity in the estimation of $z=e(t,x,w)$. More specifically, 
consider the ambiguity in the estimation of $z$ under the observation of $w$ with an error bound $\sigma$. Then the $L^p$-gain $\gamma (\sigma)$ gain equals the ratio of the ambiguity and $\sigma$. 

\subsection{Computation and example}
The input-to-output gain can be computed by solving the problem (\ref{gain}). In Section \ref{secobs}, the output function $y$ is smooth. It can be numerically approximated in a finite dimensional space, such as interpolation using a finite number of nodes. Similarly, one has to work on a finite dimensional space of $w$ to carry out the computation. So, for the purpose of computation, we only discuss the $L^p$-gain in a finite dimensional space of $w$, denoted by $\cal U$, rather than the infinite dimensional space of arbitrary integrable functions. The space ${\cal U}$ can be defined by the frequency bandwidth, or the order of polynomials, or some other spaces used for the approximation of the input. Then, (\ref{gain}) is reformulated as follows
 \EQ
\label{gain1}
\gamma_{\cal U} (\sigma)=\ds\max_{\begin{array}{c}w\in {\cal U},||w||_{L^p}\leq \sigma\\ \mu_{min}\leq \mu\leq \mu_{max}\end{array}} \Fr{||e(t,x,w,\mu)-e(t,x^\ast,0,\mu)||_{L^p}}{\sigma}
\EE
More specifically, given a positive number $\sigma >0$ define 
$$J(x(\cdot), w(\cdot),\mu)=|| e(t, x(\cdot),  w(\cdot),\mu)-e(t,x^\ast(\cdot), 0, \mu)||_{L^p}$$
Then the following dynamic optimization determines the $L^p$-gain over the space ${\cal U}$. \\
\\
{\bf Dynamic optimization for $L^p$-gain}
\EQ
\label{problemgain}
\rho=\displaystyle{\max_{( x, w, \mu)}} J\\
\hspace{-0.3in} \mbox{subject to}\\
 w(t)\in {\cal U}, \; || w(t)||_{L^p} \leq \sigma\\
\dot { x}= f(t, x,w, \mu), \\
 x(t_0)=x_0\\
 \mu_{min}\leq \mu\leq \mu_{max}
\EE
The $L^p$-gain from $w\in {\cal U}$ to $z$ is $\gamma_{\cal U}(\sigma)=\Fr{\rho}{\sigma}$. $\diamond$\\

\begin{figure}
	\begin{center}
		\includegraphics[width=2.5in]{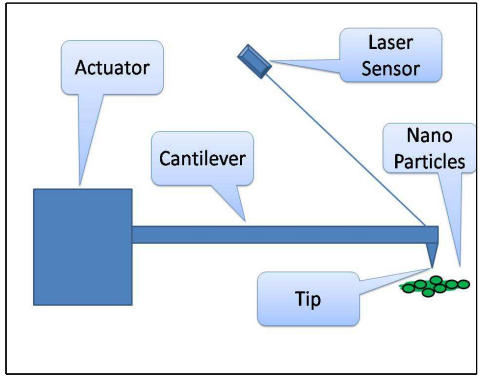} 
		\caption{Atomic force microscope}
		\label{fig:nanomanipulator}
		\end{center}
\end{figure}

{\it Example} ($L^p$-gain in the presence of system uncertainty). Atomic force microscope (AFM) invented two decades ago is used to probe surfaces at the atomic level with good accuracy. This type of equipment is also used as nano-manipulation tools to handle particles at nano-scale \cite{ashhab,basso,delnavaz}. Illustrated in Figure \ref{fig:nanomanipulator}, the system consists of a microcantilever with a sharp tip at one end. The vibration of the cantilever is measured by an optical sensor. The topographic images of surfaces can be taken by measuring the cantilever's dynamic behavior which is determined by the interacting force of the tip with the sample. 

The dynamics of the vibrating tip can be modeled as a second order system \cite{delnavaz}
\EQ
\label{nanomodel}
\dot x_1=x_2\\
\dot x_2 = -\omega^2x_1-2\xi \omega x_2+h(x_1,\delta) +u(t)+w(t)\\
z=x_1
\EE
where $x_1$ is the position of the cantilever tip at the scale of nanometers, $x_2$ is its velocity, $\omega$ is the natural frequency of the cantilever, and $\xi$ is the damping coefficient. In this system, $u(t)$ is the control input, $w(t)$ is the actuator disturbance which is unknown. The function $h(x_1,\delta)$ is the tip-sample interaction force in which $\delta$ is the separation between the equilibrium of $x_1$ and the sample surface. It is a system uncertainty. We adopt the following model for $h(x_1,\delta)$ \cite{basso}.  
$$h(x_1,\delta)=-\Fr{\alpha_1}{(\delta+x_1)^2} + \Fr{\alpha_2}{(\delta+x_1)^8}$$
At nano-scale, system uncertainly and performance robustness are critical issues in control design because a seemingly small noise or uncertainty may have significant impact on the performance. In the following, we assume that the value of $\delta$ and $w(t)$ are unknown. The goal is to compute the $L^2$-gain from the actuator disturbance $w$ to the performance $z=x_1$ for $\delta$ in the entire interval $[\delta_{min}, \delta_{max}]$.  

\begin{figure}
	\begin{center}
		\includegraphics[width=2.5in]{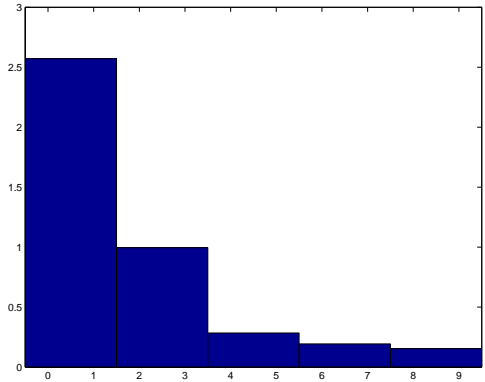} 
		\caption{$L^2$-gain}
		\label{fig:nanohist}
		\end{center}
\end{figure}

The simulations are based on the following set of parameter value
\EQQ
\omega = 1.0, & \xi = 0.02, & \alpha_1=0.1481, & \alpha_2=3.6\times 10^{-6}
\EEE
The nominal control input is $u(t)=1$ and the time interval is $[0, 7]$, which is long enough to cover one period of oscillation. The bounds of $\delta$ are
\EQQ
\delta_{min} = 0.8, & \delta_{max}=1.2
\EEE
We use $L^2$-norm as the metric for the actuator disturbance force $w$. Let $w$ be an arbitrary function in a two-frequency space ${\cal W}_{k_1,k_2}$ defined as follows
\EQQ
w=\ds\sum_{i=1}^2 \left( A_i \cos(\Fr{2\pi k_i}{t_f-t_0}t) + B_i \sin(\Fr{2\pi k_i}{t_f-t_0}t)\right)
\EEE  
Let $\sigma = 0.03$. Then the $L^2$-gain is computed by solving the dynamic optimization defined in (\ref{problemgain}). More specifically, 
\EQQ
\rho=\displaystyle{\max_{( x, u, \delta)}} || x_1(t) - x_1^\ast(t)||_{L^2}\\
\hspace{-0.3in} \mbox{subject to}\\
||w(t)||_{L^2} \leq \sigma, \; w(t)\in {\cal W}_{k_1,k_2}\\
\dot { x}= f(t, x, u, \delta), \\
\hat x(t_0)=x_0\\
\delta \in [\delta_{min}, \delta_{max}]
\EEE
The $L^2$-gain equals $\Fr{\rho}{\sigma}$. It is computed for spaces with various frequencies, ${\cal W}_{0,1}$, ${\cal W}_{2,3}$,
${\cal W}_{4,5}$, ${\cal W}_{6,7}$, and ${\cal W}_{8,9}$. The result is shown in Figure \ref{fig:nanohist}. The $L^2$-gain for frequencies $0$ and $1$ is $2.5707$. When the frequencies are increased, the gain decreases. $\diamond$\\

\subsection{An alternative algorithm for $L^2$-gain}
Solving (\ref{problemgain}) becomes increasingly difficult for high frequencies. The reason is that, for the computational purpose, the problem of dynamic optimization is always discretized at finite number of nodes in time. For higher frequencies, the number of nodes must be increased. As a result, the dimension of optimization variables is increased as well. Developing efficient methods of computation for inputs of high frequencies requires further research. 
However, in the case of $L^2$-gain, there exists an alternate approach without solving dynamic optimization. Inspired by Remark 6 and the observability gramian in \cite{krener}, the $L^2$-gain can be approximated using the following matrix. 

Suppose the space of input, ${\cal U}$, is finite dimensional with a basis, $w_1$, $w_2$, $\cdots$, $w_m$. Let $\sigma >0$ be a constant number. For the input $\pm\sigma w_i$, the trajectory of 
\EQQ
\dot x=f(t,x,\pm\sigma w_i)\\
x(t_0)=x_0
\EEE
is denoted by $x^{i\pm}(t)$. Define
\EQ
\label{dltz}
2\Delta z^i=e(t, x^{i+}(t),\sigma w_i(t))-e(t, x^{i-}(t),-\sigma w_i(t))
\EE
Now, define 
\EQ
\label{Gmatrix}
G^w_{ij}=<w_i, w_j>=\Fr{1}{t_f-t_0} \ds\int_{t_0}^{t_f} w_i(t)^Tw_j(t)dt\\
G^z_{ij}=<\Delta z^i,\Delta z^j>=\Fr{1}{t_f-t_0} \ds\int_{t_0}^{t_f} \Delta z^i(t)^T\Delta z^j(t)dt
\EE
Denote the matrices $G^w=(G^w_{ij})_{i,j=1}^n$ and $G^z=(G^z_{ij})_{i,j=1}^n$. Given any 
$$w=\ds\sum_{i=1}^m a_iw_i$$
satisfying 
$$<w,w>=\sigma^2$$
Then $\Delta z$ generated by $\pm w$ is approximately
$$2\Delta z=e(t, x^{+}(t),w(t))-e(t, x^{-}(t), -w(t))\approx 2\ds\sum_{i=1}^m a_i\Delta z^i$$
Therefore, 
\bea
||w||_{L^2}^2 &=& \MT a_1 \cdots a_m\EM G^w\MT a_1 & \cdots & a_m\EM^T\\
||\Delta z||_{L^2}^2 &\approx& \MT a_1 \cdots a_m\EM G^z\MT a_1 & \cdots & a_m\EM^T
\eea
Therefore, the $L^p$-gain square is approximately the solution of the following optimization
\EQQ
\Fr{1}{\sigma^2} \max_a a^T G^z a\\
\hspace{-0.3in} \mbox{ subject to}\\
a^T G^w a=\sigma^2
\EEE
where $a=\MT a_1&a_2&\cdots & a_m\EM^T$. A necessary condition for the optimal solution is
$$G^z a=\lambda G^w a$$
for some scalar $\lambda$. At this point, 
\bea 
a^TG^za\==\lambda a^TG^wa\\
\== \lambda \sigma^2
\eea
Therefore, 
\bea
\gamma_{\cal U}(\sigma)^2\==\ds\max_{||w||_{L^2}=\sigma^2} \Fr{|| \Delta z||_{L^2}^2}{||w||_{L^2}^2}\\
&\approx& \Fr{a^TG^za}{\sigma^2}\\
\== \lambda
\eea 
On the other hand, $\lambda$ is an eigenvalue of $(G^w)^{-1}G^z$. So, the $L^2$-gain is approximately the square root of the largest eigenvalue.

To summarize, given a system 
\EQQ
\dot x=f(t,x,w)\\
x(t_0)=x_0
\EEE
and a space of input functions  ${\cal U}$ with basis $w_1$, $w_2$, $\cdots$, $w_m$. Given $\sigma >0$, compute $\Delta z^i$ in  (\ref{dltz}). Compute the matrices $G^w$ and $G^z$ in (\ref{Gmatrix}). Then the $L^2$-gain is approximately $\sqrt{\lambda_{max}}$, where $\lambda_{max}$ is the largest eigenvalue of $(G^w)^{-1}G^z$. 
\\

\begin{figure}
	\begin{center}
		\includegraphics[width=2.5in]{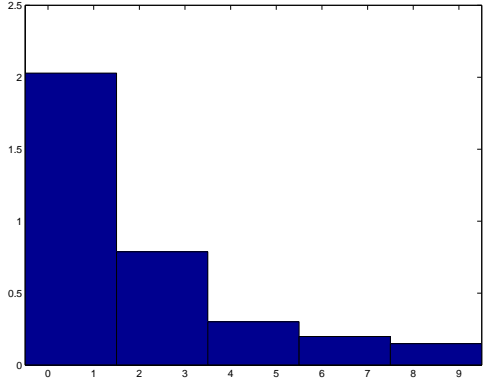} 
		\caption{Approximate $L^2$-gain}
		\label{fig:nanohistgrammian}
		\end{center}
\end{figure}

\noindent{\it Example}. 
Consider the model of AFM defined in (\ref{nanomodel}). Assume that the value of $\delta=1.0$ is known. The other parameters are the same as in the previous example. The approximate $L^2$-gain is computed using the matrix approach. The result is shown in Figure \ref{fig:nanohistgrammian}. This method is straightforward in computation because no optimization is required. However, the approximation does not take into full account the nonlinear dynamics. Comparing to the gain using (\ref{problemgain}), the errors of the gain computed using covariance matrix is around $18 \sim 20\%$ in ${\cal W}_{0,1}$, ${\cal W}_{2,3}$, and ${\cal W}_{4,5}$.

\section{Reachability}
\label{seccon}
Dynamic optimization can be applied to quantitatively measure reachability. Consider a control system 
\bean
\label{sys2}
\dot x=f(t,x,u)
\eean
where $x \in \Re^{n_x}$ and $ u\in \Re^{n_u}$. Suppose $||x||_{X}$ is a norm in $\Re^{n_x}$. We suppose that the state and control are subject to constraint
$$(x(\cdot), u(\cdot))\in {\cal C}$$ 

\begin{definition}
\label{defcon}
Given $x_0$ and $x_1$ in $\Re^{n_x}$.  Define
\EQ
\label{rhox0}
\rho_c(x_0, x_1)^2=\ds\min_{(x, u)}||x(t_1)-x_1||^2_{X}\\
\hspace{-0.3in} \mbox{subject to}\\
\dot x=f(x,u)\\
x(t_0)=x_0\\
(x(\cdot),u(\cdot)) \in \cal C
\EE
The number $\rho_c(x_0, x_1)$ is called the ambiguity in control.  

Let $D_0, D_1 \subset \Re^{n_x}$ be subsets in state space. The ambiguity in control over the region $\bar D_0\times \bar D_1$ is defined by the following max-min problem.
$$\begin{array}{lll}
\rho_c=\displaystyle{\max_{(x_0, x_1) \in \bar D_0\times \bar D_1}}\rho_c(x_0, x_1)
\end{array}$$
\end{definition} 

In this definition, $t_1$ is either fixed or free in a time interval $[t_0, T]$. In the following discussion, we assume $t_1$ is fixed. If $\rho_c(x_0,z_1)$ is nonzero, then the state cannot reach $x_1$ from $x_0$ by using admissible controls. The maximum value of $\rho_c(x_0,z_1)$ over $\bar D_0\times \bar D_1$ represents the worst  scenario of reachability. In some applications, the relative ambiguity
$$\Fr{\rho_c(x_0,x_1)}{||x_1||_{X}}$$ 
is used to measure the reachability.

The definition of $\rho_c$ is consistent with the classic definition of controllability for linear time-invariant systems. To be more specific, consider a linear system
\bean
\label{syslinear}
\dot x=Ax+Bu
\eean
let $u(\cdot)\in {\cal C}$ be the space of continuous functions from $[t_0, t_1]$ to $\Re^m$; and let $D_0=D_1=\Re^n$. If $(A,B)$ is controllable, i.e. 
$$rank(\MT B&AB&A^2B&\cdots &A^{n-1}B\EM=n$$
then for any $x_0$ and $x_1$, there always exists a control input so that $x(t)$ with $x(0)=x_0$ reaches $x_1$ at $t=t_1$. Therefore, $\rho_c(x_0,x_1)$ is always zero for arbitrary $(x_0,x_1)$; and $\rho_c = 0$. On the other hand, if $(A,B)$ is uncontrollable, then under a change of coordinates an uncontrollable subsystem can be decoupled from the controllable part of the system. In the uncontrollable subsystem, the states cannot be driven to close to each other by control inputs. Therefore, $\rho_c(x_0,x_1)$ is unbounded for arbitrary states in the uncontrollable subspace. This implies $\rho_c=\infty$. To summarize,
$$\rho_c=\left\{\begin{array}{lll}
0 & \mbox{if } (A,B)\mbox{ is controllable }\\
 \infty & \mbox{if } (A,B)\mbox{ is uncontrollable }
 \end{array}\right.
$$

Control has a cost. For weakly reachable systems, it takes relatively large control energy to reach a terminal state. The cost in reachability can be measured by the following quantity. Denote $||u(\cdot)||_{\cal U}$ and $||x||_X$ the metrics of the control input and the state, respective. 

\begin{definition}
Given initial and final states, $x_0$ and $x_1$, define
\EQ
\label{rhox0b}
W(x_0, x_1)=\ds\min_{(x, u)}\ds\lim_{\psi\rightarrow \infty}\left(||u(t)-u^\ast(t)||_{\cal U}+\psi ||x(t_1)-x_1||_X\right)\\
\hspace{-0.3in} \mbox{subject to}\\
\dot x=f(x,u)\\
x(t_0)=x_0,\\
(x(\cdot),u(\cdot)) \in \cal C
\EE
\end{definition}

This definition has the following property. If $x_1$ can be reached from $x_0$, then $W(x_0,x_1)$ equals the minimum control
\EQQ
W(x_0, x_1)=\ds\min_{(x, u)}||u(t)-u^\ast(t)||_{\cal U}\\
\hspace{-0.3in} \mbox{subject to}\\
\dot x=f(x,u)\\
x(t_0)=x_0,\;\; x(t_1)=x_1\\
(x(\cdot),u(\cdot)) \in \cal C
\EEE
If $x(t)$ cannot reach $x_1$ using admissible control, then $W(x_0, x_1)=\infty$. A large value of $W(x_0,x_1)$ implies higher control cost, thus weak reachability.\\

\noindent \underline{\it Remark 7}. 
Suppose the system is linear. Suppose $D_0=\{ 0\}$ and $D_1=\{ x| \;\;||x||<\epsilon\}$ for some small $\epsilon >0$. Let $W$ be the maximum cost in reachability under $L^2$-norm, 
$$W=\ds\max_{x_1\in \bar D_1}W(0,x_1)$$
Then $(W/\epsilon)^2$ equals the reciprocal of the smallest eigenvalue of the controllability gramian. 

To justify Remark 7, consider 
$$\dot x=Ax+Bu$$
If $(A,B)$ is uncontrollable, we know $W=\infty$. We also know that the smallest eigenvalue of the controllability gramian is zero. Therefore the claim holds true. Now suppose $(A,B)$ is controllable. Define
$$|| u||_{\cal U}^2=\int_{t_0}^{t_1} ||u(t)||^2dt$$
The control cost to reach $x_1$ is defined by 
\EQQ
W(0, x_1)=\ds\min_{(x, u)}||u(t)||_{\cal U}\\
\hspace{-0.3in} \mbox{subject to}\\
\dot x=Ax+Bu\\
x(t_0)=0, x(t_1)=x_1
\EEE
Let $P$ be the controllability gramian, then it is known \cite{zabczyk} that the optimal cost satisfies
\bea
&&\ds\int_{t_0}^{t_1} ||u(t)||^2dt\\
\== x_1^T P^{-1}x_1
\eea
If $\sigma_{min}$ is the smallest eigenvalue of $P$, then
\bea
\label{eqcomp}
(W/\epsilon)^2&\==&\Fr{1}{\epsilon^2}\max_{||x_1||\leq \epsilon} x_1^T P^{-1}x_1\\
&\==& \Fr{1}{\sigma_{\min}} 
\eea
$\diamond$

\subsection{Example} 
In the following, we compute the ambiguity in the control of a heat equation with Neumann boundary control.  
\bean
\label{pdemodel}
&&\Fr{\partial w(r,t)}{\partial t} -\kappa \Fr{\partial^2 w(r,t)}{\partial r^2}=0\\
&&w(r,0)=0, \;\;\; 0 \leq r\leq 2\pi\\
&&w(0, t)=0, \;\;\; 0\leq t\leq t_f\\
&&w_{rt}(r,t)|_{r=2\pi}=u(t)
\eean
where $w(r,t)\in \Re$ is the state of the system, $r\in \Re$ is the space variable, and $t$ is time. The control input is $u$. For a thermal problem, $u$ represents the rate of the heat flux $w_r$. The initial state is assumed to be zero. 

For the purpose of computation, we discretize the problem at equally spaced nodes, 
$$0=r_0< r_1<r_2<\cdots < r_N=2\pi$$
Define 
$$x_1(t)=w(r_1,t), \; x_2(t)=w(r_2, t), \; \cdots, x_N(t)=w(r_N,t)$$ 
Using central difference in space, (\ref{pdemodel}) is approximated by the following control system defined by ODEs.
\bean
\label{pdefinited}
\dot x_1 \== \kappa \Fr{x_2-2x_1}{\Delta r^2}\\
\dot x_2\== \kappa \Fr{x_{1}+x_{3}-2x_2}{\Delta r^2},\\
&\vdots&\\
\dot x_i\== \kappa \Fr{x_{i-1}+x_{i+1}-2x_i}{\Delta r^2},\\
&\vdots\\
\dot x_{N-1}\==\kappa \Fr{x_{N-2}+x_{N}-2x_{N-1}}{\Delta r^2}\\
\dot x_N\==v\\
v\==\dot x_{N-1}+\Delta r u
\eean

\begin{figure}
	\begin{center}
		\includegraphics[width=3.0in]{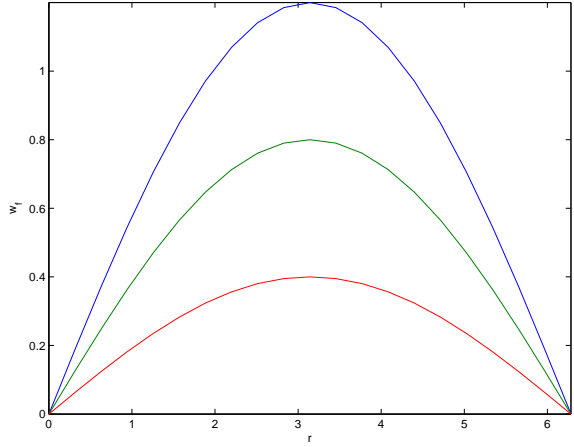} 
		\caption{The target state $w_f(r)$}
		\end{center}
	\label{figarch}
\end{figure}

We understand that more sophisticated algorithms of solving the heat equation exit. We adopt this central difference method for simplicity in the illustration of control ambiguity. System (\ref{pdefinited}) is linear and controllable. So, it is theoretically a reachable system. However, in reality the maximum temperature cannot exceed safety margin. Under such constraint, a controllable linear system may not be reachable due to overshot. Let $t_f=150$ and $\kappa=0.14$. Suppose the target states are the following arches. A few of them is shown in Figure \ref{figarch}.
$$w_f(r)=w(r,t_f)=A\sin (r/2), \;\;\; 0\leq r\leq 2\pi$$
The goal is to compute the ambiguity in control from $w(r,0)=0$ to $w_f(r)$ with the magnitude of $0\leq A\leq 1.2$ subject to the constraint 
$$w(r,t)\leq 2$$ 
The norm in the finite dimensional state space is defined by 
$$||x||^2=\Delta r\sum_{i=1}^N  x_i^2$$
which approximates the $L^2$-norm in $C[0, 2\pi]$. To compute the control ambiguity in the range of $0\leq A\leq 1.2$, we consider the following nodes in the magnitude of the target arch
$$A=0, \;0.2, \;0.4, \;0.6, \;0.8, \;1.0, \;1.2$$ 

\begin{figure}
	\begin{center}
		\includegraphics[width=3.0in]{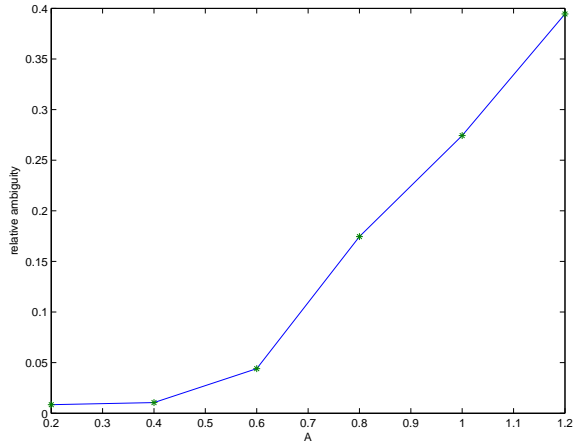} 
		\caption{Magnitude of target state vs. relative control ambiguity}
		\end{center}
	\label{fig:constrained}
\end{figure}

The corresponding dynamic optimization problem for the control ambiguity is defined by 
\EQ
\label{pderho}
\rho_c(0, x_f)^2=\ds\min_{(x, u)}\Delta r\sum_{i=1}^N  (x_i-w_f(r_i))^2\\
\hspace{-0.3in} \mbox{subject to}\\
\dot x=f(x,v)\\
x_i(t)\leq 2,\;\;\; i=1,2,\cdots,N\\
x(0)=0
\EE
where $f(x,v)$ is the dynamics defined in (\ref{pdefinited}) with the input $v$, and  
$$x_f=\MT w_f(r_1)&w_f(r_2)&\cdots &w_f(r_N)\EM$$
In the simulation, $N$ is selected to be $N=31$. Problem (\ref{pderho}) is solved using Pseudospectral method at Legendre-Gauss-Lebato (LGL) nodes \cite{EKR1995,FR, KGR}. We use $15$ LGL nodes in this example. 
Through computation, it is found that the system becomes increasingly unreachable due to the constraint when the value of $A$ is bigger than $0.4$. The relative ambiguity in control is shown in Figure \ref{fig:constrained}. When the magnitude of the target state is $1.2$, the relative ambiguity shows that the closest state that the system can reach has almost a $40\%$ relative error. Therefore, the system is practically unreachable if the state is required to be bounded. 

\section{Conclusion}
It is shown by numerous examples and definitions that computational dynamic optimization is a promising tool of quantitatively analyzing control system properties. Using computational approaches, the concepts studied in this paper, including the ambiguity in estimation and control, input-to-output gain, and the cost in reachability, are applicable to a wide spectrum of applications. In addition, these concepts are defined and applied in a way so that one can take advantage of user knowledge or take into account system constraints. As a result, the properties of control systems are not only verified, but also measured quantitatively. While these concepts can be applied to a wide spectrum of problems, some specific applications exemplified in this study include: strongly observable (detectable) or weakly observable (detectable) systems; improving observability by employing user knowledge; partial observability of networked complex systems; $L^2$-gain of nonlinear control systems; reachability in the presence of state constraints; and boundary control of partial differential equations.  

Similar to many nonlinear optimization problems, a main drawback of the approach is that a global optimization is, in general, not guaranteed for nonlinear systems. In addition, the problem of computational accuracy also poses many questions remain to be answered.

\end{document}